\documentclass{amsart}

\newcommand\thisuni{Maynooth University}%for published version
\newcommand\acolleague{John Murray }%for pub version
\newcommand\afriend{Des MacHale}%for pub version
\newcommand\collaborators{Tirthankar Bhattacharyya and collaborators }%pub v
\newcommand\ba{\begin{aligned}}
\newcommand\ea{\end{aligned}}
\newcommand\be{\begin{equation}}
\newcommand\ee{\end{equation}}
\newcommand\C{\mathbb{C}}
\newcommand\GL{\textup{GL}}
\newcommand\ignore[1]{}
\newtheorem{question}{Question}

\newcommand{\R}{\mathbb{R}}
\newcommand{\half}{{\scriptstyle\frac12}}

\title{Forms of nice questions}
\author{Anthony G. O'Farrell} 
\date{\today}
\address{Mathematics and Statistics, Maynooth University, Co Kildare, Ireland W23 HW31}
\email{anthony.ofarrell@mu.ie}

\RequirePackage{url}
\ignore{
	\authorbio{
Anthony O'Farrell 
studied at University College Dublin and Brown University, and worked at UCLA
before taking the chair of Mathematics in Maynooth, Ireland, where he served for
37 years and is now Professor Emeritus.  His published work is mainly in analysis, but includes a monograph about reversibility in dynamics and group theory
and a Mathematical Olympiad manual. 
More at \url{https://www.logicpress.ie/aof}.
}
}%end ignore

\begin{document}

\begin{abstract}
You can invent striking and challenging problems with unique
solution by building some symmetry into functional equations.
Some are suitable for high school; others could generate
college-level projects involving computer algebra. The
problems are functional equations with group actions in the
background.  Interesting examples arise even from small finite
groups. Whether a given problem ``works" with a given choice of 
constant coefficients depends on whether a related multilinear form is nonzero.
These forms are essentially the classical group determinants
studied by Frobenius in the nineteenth century.
%86 words
\end{abstract}

\maketitle

\section{Nice questions}\footnote{The version of record of this article
will appear in the College Mathematics Journal, with DOI 10.1080/07468342.2025.2485013.
It will be accessible on open access from the Taylor and Francis website.}
At \thisuni, we run ``Mathematical enrichment sessions,'' aimed at
encouraging participation in the Irish and International 
Mathematical Olympiads.
Recently, a colleague called for problems that might
be useful for a practice contest.  
I proposed the following question.

\begin{question}\label{Q:1}
Suppose $f(\tan\theta)+ 2f(\cot\theta) = \cos(2\theta)$ for $0<\theta<\pi/2$.
Find $f(2024)$.
\end{question}

Substituting $t$ for $\tan\theta$, this becomes
\be\label{E:1-1}
 f(t) + 2 f\left(\frac1t\right) = \frac{1-t^2}{1+t^2}.
\ee
Replacing $t$ by $1/t$ gives
\be\label{E:1-2}
 f\left(\frac1t\right) + 2f(t) = \frac{t^2-1}{1+t^2}.
\ee
Subtracting \eqref{E:1-1} from twice \eqref{E:1-2} we get
a formula for $f(t)$,
$$ f(t) = \frac{t^2-1}{t^2+1}, $$
which gives $f(2024)$ as $4096575/4096577$.

What happened here is that a single linear equation for an unknown
function generated a linear system that had a unique solution. 
This came about because of a symmetry in the equation:
the function $\tau:t\mapsto 1/t$ is an \emph{involution} on the 
set $(0,+\infty)$ of positive real numbers,
i.e. it is a function that inverts itself.   
 
\subsection{} 
There are plenty of involutions $\tau$ on all or part of the real line
(or the complex plane). For instance the reflections $x\mapsto -x$
and $x\mapsto 2-x$ are involutions on $\R$, and
$x\mapsto -x/(x+1)$ is an involution on the line or plane with 
$x=-1$ deleted.  So precisely the same idea allows us to propose any problem
in any of forms
\begin{eqnarray}
	a f(x) + b f(-x) &= F(x),\label{E:2-1}\\
	a f(x) + bf(2-x) &= F(x),\label{E:2-2}\\
	a f(x) + b f\left(\frac{-x}{x+1}\right) &=F(x),\label{E:2-3}
\end{eqnarray}
where we specify $a,b$ and $F(x)$ and ask for the unknown $f(x)$.
There is going to be a unique solution except when $a^2=b^2$.
One may use variable coefficients $a(x)$ and $b(x)$, provided
one avoids points where $a(x)a(\tau(x))=b(x)b(\tau(x))$.

\subsection{}
A similar situation
also surfaced in recent work with
\collaborators \cite{BORK}.
We studied the equation
\begin{equation}\label{E:3-1}
 f(z) + z\,\overline{f(z)} = F(z),
 \end{equation}
for $z$ in the open unit disc. The action of the involution
$z\mapsto \bar z$, complex conjugation, this time on the image
instead of the argument, gives the second equation
$$ \overline{f(z)} + \bar z f(z) = \overline{F(z)}. $$
Since the determinant $1-|z|^2$ does not vanish on the disc,
it follows that the original equation has unique solution $f(z)$
for each given complex-valued function $F(z)$, given by
$$ f(z) = \frac{F(z)-z\overline{F(z)}}{1-|z|^2}.  $$

\subsection{} 
The following case 
 does not involve an involution. 
Let $\omega=\frac12+\frac{\sqrt3}2i$, one of the complex cube roots 
of unity, and consider the following question.

\begin{question}\label{Q:3}
	Suppose $f:\C\to\R$ and
\begin{equation}\label{E:4-1}
f(z) - f(\omega z) + f(\omega^2 z) = z^2.
\end{equation}
	Find $f(10)$.
\end{question}

Replacing $z$ in turn by $\omega z$ and $\omega^2 z$ gives 
two more equations, and together with \eqref{E:4-1} we have a
system of three linear equations connecting the unknowns
$f(z), f(\omega z), f(\omega^2 z)$. We can eliminate the latter
two in the usual way, and obtain
$$
f(z) = -\half \omega z^2,  
$$
so that $f(10)=-50\omega= 25(1+\sqrt{-3})$.
\ignore{A devious person could disguise the r\^ole of $\omega^2 z$
in this kind of problem by omitting the third term on the left.
The problem can still be solved, in the same way.
}

\section{Group actions}
There are many more questions one can pose along these lines.
All these questions involve some kind of symmetry, and the mathematics
of symmetry is group theory. To describe the general framework
that prompted these questions, we recall the terminology of
group actions.

An \emph{action} of a  group $G$ on a set $X$
is a homomorphism from $G$ into the group of permutations of $X$
(under composition).
In other words, an action $\alpha$ associates to each element $g\in G$
a bijection $\alpha(g):X\to X$, and for any two elements $g_1,g_2\in G$,
we have
$$ \alpha(g_1g_2) = \alpha(g_1)\circ \alpha(g_2), $$
i.e., for all $x\in X$,
$$ \alpha(g_1g_2)(x) = \alpha(g_1) \left( \alpha(g_2) (x) \right) 
.$$
The action is \emph{faithful} if $\alpha$ is injective. So
a faithful action just identifies an isomorphic copy of $G$ inside
the permutation group on $X$. When we have such an action, we simplify the 
notation by identifying $G$ with its image, writing $g(x)$ instead
of $(\alpha g)(x)$. When $X$ has a vector space structure and $\alpha(g)$
is linear, we simplify further by writing $gx$ instead of $g(x)$.

The functional equations we are considering are associated to 
actions of finite groups on the domain or the image of
the unknown function, or on both domain and image.

\subsection{}
For example, take the group $C_2$ of order $2$,
generated by a single
element $\tau$ having $\tau^2=1$. A faithful action of $C_2$ on $\R$
identifies $\tau$ with some involution on $\R$.  
Each such action gives us a linear map
$(a,b) \to L(a,b)$ from $\R^2$ to operators on real-valued functions
of a real variable, where 
$$ L(a,b) f := a\cdot f + b\cdot  (f\circ\tau). $$ 
The examples  in equations \eqref{E:1-1},\eqref{E:2-1}, \eqref{E:2-2}, and \eqref{E:2-3}, 
work because $L$ is bijective if
$a^2\not=b^2$. The single equation
\begin{equation*}%\label{E:5-1}
a f(x) + b f(\tau(x)) = F(x)
\end{equation*}
(with given $F$ and unknown $f$) 
is equivalent to the second equation
\begin{equation*}%\label{E:5-2}
a f(\tau(x)) + b f(x) = F(\tau(x)). 
\end{equation*}
Provided $a^2\not= b^2$, one can combine the two to
eliminate $f(\tau(x))$,
and write $f(x)$ as a linear combination of $F(x)$ and $F(\tau(x))$.

In Equation \eqref{E:2-3}, the map $x\mapsto -x/(x+1)$
has a singularity at $x=-1$, so the action is on $\R\setminus\{-1\}$,
instead of on the whole real line.

\subsection{} 
The group of order 3 is denoted $C_3$. 
The problem in 
Question \ref{Q:3} is derived from a faithful
action of $C_3$ on the complex plane. 
The rotation $\sigma:z\mapsto \omega z$ generates
an isomorphic copy of $C_3$ in the group of bijections of the complex
plane, so the single equation
$$ a f + b f\circ\sigma + c f\circ \sigma\circ\sigma = F $$
is equivalent to the $3\times 3$ system
\begin{equation*}\nonumber\ba
	a f + b f\circ\sigma + c f\circ \sigma\circ\sigma &= F\\ 
	c f + a f\circ\sigma + b f\circ \sigma\circ\sigma &= F\circ\sigma\\ 
	b f + c f\circ\sigma + a f\circ \sigma\circ\sigma &= F\circ\sigma\circ\sigma. 
	\ea
\end{equation*}
This has a unique solution for every complex-valued function $F$
of a complex variable, provided the determinant of this linear
system, the circulant
$$
\left|\begin{matrix} a&b&c \\ c&a&b \\ b&c&a \end{matrix}\right|
= a^3+b^3+c^3-3abc $$
is nonzero. 

There are many nice actions of $C_3$ on the complex plane, generated by
particular linear fractional transformations 
(to be precise, one usually has to puncture the plane, removing two or three
points).
An example is generated by the map 
$\displaystyle\sigma(z)=\frac{-1}{z+1}$. 
%or $-(2z+7)/(z+3)$...
%Some readers will recognise elements of PSL(2,\Z)
This example maps reals to reals
(and even rationals to rationals), so can be used in elementary classes 
without involving complex numbers.

\medskip
More generally, an action of any cyclic group $C_n$ of order $n$
gives problems in which a single equation generates an $n\times n$
linear system having a circulant for its determinant.
For instance $z\mapsto iz$ generates an action of $C_4$ on the plane.

\subsection{} 
The example in Equation \eqref{E:3-1}
involves a group
acting on image of the unknown function $f$.
The complex conjugation map $\kappa:z\mapsto \bar z$
generates an action of $C_2$ on the plane, and the equation 
is of the form 
\begin{equation}\label{E:6-1}
	a\cdot f + b\cdot (\kappa\circ f) = F, 
\end{equation}
where $a,b,F$ are given complex-valued functions on
some domain --- the unit disk in this case --- and $f$ is unknown.
But $\kappa$ is not just any involution on the plane;
it is a ring automorphism of the complex number ring $(\C,+,\cdot)$.
As a result,
$$\kappa\circ(f\cdot g) = (\kappa\circ f)\cdot (\kappa\circ g)$$ and
$$\kappa\circ(f + g) = (\kappa\circ f) + (\kappa\circ g),$$ 
i.e., the composition
$f\mapsto \kappa\circ f$ is an involutive ring automorphism
of the ring of all complex-valued functions on the domain.
Applying the composition to both sides, Equation \eqref{E:6-1} is
equivalent to the $2\times2$ system
\begin{equation*}\nonumber\ba
	a\cdot f + b \cdot (\kappa\circ f) &= F\\ 
	(\kappa\circ b)\cdot  f + (\kappa\circ a)\cdot (\kappa\circ f) 
	&= \kappa\circ F, 
	\ea
\end{equation*}
so we get unique solutions provided
the determinant 
$$ a\cdot (\kappa\circ a) - b\cdot (\kappa\circ b) $$
does not vanish on the domain in question.

In fact, complex conjugation is the only continuous nonidentity automorphism 
of the ring of complex numbers \cite{Yale}, so for complex-valued functions
this method applies only to cases where $|a|\not=|b|$ on the domain
in question.  However, more possibilities arise when the functions
take other kinds of values. 

\subsection{The general case.}\label{SS:general}
Combining all these ideas, we arrive at the following framework.

Consider functions $f:X\to A$ where a finite group $G$ acts on the set $X$,
and a finite group $H$ acts on a field $F$
and on an associative algebra 
$A$ over $F$ 
(\cite[Chapter 6]{Herstein}, \cite[Chapter 7]{Jacobson}, 
\cite[Chapter 8]{Robinson}, \cite[Chapter 6]{Rotman}
or \cite[Chapter 8]{Fraleigh})  
so that for each $h\in H$ and for all $a\in F$ and all $y\in A$,
$$ h(a y) = (ha) hy, 
$$
and for all $y_1,y_2\in A$,
$$ h(y_1+y_2) = hy_1 + hy_2.$$ 
Then a single linear functional equation, valid for all $x\in X$,
\begin{equation}\label{E:general}
\sum_{g\in G, h\in H} a_{g,h} h( f( g(x) ) ) = F(x),  
\end{equation}
(where the coefficients $a_{g,h}\in F$ are given for $g\in G$ and $h\in H$)
will have unique solution $f$ for each given $F:X\to A$,
whenever
the form in $|G|\cdot|H|$ variables of degree $|G|\cdot|H|$
given by the determinant $\Delta$ of the coefficients of the system
of functional equations
\begin{equation}\label{E:general-system}
\sum_{g\in G, h\in H} (ra_{g,h}) (rh( f( g(k(x))) ) ) = rF(k(x))
, 
\end{equation}
(where $k$ ranges over all elements of $G$ 
and $r$ ranges over all elements of $H$)
is nonzero.
This is a homogeneous form with coefficients in $F$,
involving $|G|\cdot|H|$ variables.
(The forms that arise in this way are obliquely referenced
in the title of the present  article,  which is intended as
a gentle pun.)

We give a few examples of these forms.

\subsection{} 
The symmetric group on three symbols is denoted $S_3$.
The group
$G=S_3$ acting on the domain 
and the trivial group $H=(1)$ acting on the image results
in the determinant 
$$\left|
\begin{matrix}
        a&b&c&d&e&f\\
        b&a&f&e&d&c\\
        c&e&a&f&b&d\\
        d&f&e&a&c&b\\
        f&d&b&c&a&e\\
        e&c&d&b&f&a
\end{matrix}
\right|.
$$
This evaluates to\footnote{I have provided some detail on
the calculations in this paper in the document
\textit{Supplement to Forms of Nice Questions}
which may be downloaded from my publications page
at \url{https://www.logicpress.ie/aof/publications.html}.
See item 21 in the list of Expository Papers.}
\begin{equation*}%\label{E:S3-1}
	\ba
(a^2 - b^2 + bc - c^2 + bd + cd - d^2 - ae + e^2 - af - ef + f^2)^2&\\
\cdot (a + b + c + d + e + f)(a - b - c - d + e + f).&
\ea\end{equation*}
\ignore{
Sage:
var ('a,b,c,d,e,f');

A4 = matrix(SR,6,6,[a,b,c,d,e,f,b,a,f,e,d,c,c,e,a,f,b,d,d,f,e,a,c,b,f,d,b,c,a,e,e,c,d,b,f,a]);

A4; A4.determinant().expand().factor();

Result:
(a^2 - b^2 + b*c - c^2 + b*d + c*d - d^2 - a*e + e^2 - a*f - e*f + f^2)^2*(a + b + c + d + e + f)*(a - b - c - d + e + f)

}%end ignore

One example of an action of $S_3$ on the complex plane identifies
the group with the group generated by the 
linear fractional maps $1-t$ and $1/t$.
This give equations such as
\begin{equation*}%\label{E:S3-2}
\ba
h(t)+2h(1-t) + 3h\left(\frac1t\right)
+ 4h\left(\frac{1}{1-t}\right)\ &\\
+ 5h\left(\frac{t}{t-1}\right)
+ 6h\left(\frac{t-1}{t}\right)
&= F(t),
\ea
\end{equation*}
for which the determinant is $3024\not=0$.
\ignore{
	var ('a,b,c,d,e,f,V');
def find(v):
    a=1
    for b in range(2,7):
        for c in range(2,7):
            for d in range(2,7):
                for e in range(2,7):
                    for f in range(2,7):
                        V=Set([a,b,c,d,e,f])
                        if (V.cardinality()==6):

                            A4 = matrix(SR,6,6,
                            [a,b,c,d,e,f,
                             b,a,d,c,f,e,
                             c,f,a,e,d,b,
                             f,c,e,a,b,d,
                             e,d,f,b,a,c,
                             d,e,b,f,c,a]);
                            if (A4.determinant()==v):
                                print(a,b,c,d,e,f)
                                return

find(3024)
This determinant differs from that on the previous page in that d and e have been interchanged.
}%end ignore
Here, we may consider complex, or real, or just rational
variables.

One may also use coefficients from the field of rational functions,
and consider equations such as
\begin{equation}\label{E:S3-3}
	(1+t)h(t)+(1-t)h(1-t) + \frac1th\left(\frac1t\right) =  
F(t).
\end{equation}
The six(!) coefficients are now $a=1+t$, $b=1-t$, $c=1/t$
and $d=e=f=0$. They change to other rational functions
when we compose both sides
with elements of the group, but we still get a $6\times6$ linear system,
and we find that the determinant equals 
$-4$, a nonzero rational function.
\ignore{
Sage:
var('t')
A5 = matrix(SR,6,6,[1+t,1-t,1/t,0,0,0,
        t,1+1-t,0,1/(1-t),0,0,
        t,0,1+1/t,0,0,1-1/t,
        0,1-t,0,1+1/(1-t),t/(t-1),0,
        0,0,0,1/(1-t),1+t/(t-1),1-1/t,
        0,0,1/t,0,t/(t-1),1+(t-1)/t]);

print(A5);
A5.determinant().expand().factor()

gives
-4
}%end ignore
Thus Equation \eqref{E:S3-3} will have a unique rational function solution
$h(t)$ for each given rational function $F(t)$.

%\subsection{} 
An example of the case $G=C_2, H=C_2$, with group actions on
both domain and image,
is the equation
\begin{equation*}
a f(z) + bf(\bar z) +c\overline{f(z)}
+ d\overline{f(\bar z)} =F(z),
\end{equation*}
where $z$ is a complex variable.

The determinant is
$$\left|
\begin{matrix}
        a&b&c&d\\
        b&a&d&c\\
        \bar c&\bar d&\bar a &\bar b\\
        \bar d&\bar c & \bar b & \bar a
\end{matrix}
\right|
= (|a+b|^2-|c+d|^2)
(|a-b|^2-|c-d|^2). $$

\subsection{}
\bigskip
The quaternion $8$-group acts by
left-multiplication on the quaternions. 
The equation for a function $f$ of a quaternion variable
(or a variable in any space on which the quaternion
group acts) takes the following form.
\begin{equation*}\ba
a_1f(x)+ a_{-1}f(-x) + a_if(ix) + a_{-i} f(-ix)&\\
+
a_jf(jx) + a_{-j}f(-jx)
+a_kf(kx) + a_{-k}f(-kx)
&= F(x).
\ea
\end{equation*}
The associated determinant form is
$$\left|
\begin{matrix}
        b_1&b_2&b_3&b_4&b_5&b_6&b_7&b_8\\
        b_2&b_1&b_4&b_3&b_6&b_5&b_8&b_7\\
        b_4&b_3&b_1&b_2&b_7&b_8&b_6&b_5\\
        b_3&b_4&b_2&b_1&b_8&b_7&b_5&b_6\\
        b_6&b_5&b_8&b_7&b_1&b_2&b_3&b_4\\
        b_5&b_6&b_7&b_8&b_2&b_1&b_4&b_3\\
        b_8&b_7&b_6&b_5&b_3&b_4&b_1&b_2\\
        b_7&b_8&b_5&b_6&b_4&b_3&b_2&b_1
\end{matrix}
\right|,
$$
where we have relabelled 
$a_1$, $a_{-1}$, $a_i$, $a_{-i}$, $a_j$, $a_{-j}$, $a_k$, and $a_{-k}$
as $b_1$, $b_2$, $b_3$, $b_4$, $b_5$, $b_6$, $b_7$, and $b_8$, respectively.
This equals
\begin{equation*}\nonumber\ba
&(b_1^2 - 2b_1b_2 + b_2^2 + b_3^2 - 2b_3b_4 + b_4^2 + b_5^2 - 2b_5b_6 + b_6^2 + b_7^2 - 2b_7b_8 + b_8^2)\\
&\cdot(b_1^2 - 2b_1b_2 + b_2^2 - b_3^2 + 2b_3b_4 - b_4^2 - b_5^2 + 2b_5b_6 - b_6^2 + b_7^2 - 2b_7b_8 + b_8^2)\\
&\cdot(b_1 + b_2 + b_3 + b_4 + b_5 + b_6 + b_7 + b_8)\\
&\cdot(b_1 + b_2 + b_3 + b_4 - b_5 - b_6 - b_7 - b_8)\\
&\cdot(b_1 + b_2 - b_3 - b_4 + b_5 + b_6 - b_7 - b_8)\\
&\cdot(b_1 + b_2 - b_3 - b_4 - b_5 - b_6 + b_7 + b_8).
\ea\end{equation*}
or
\begin{equation*}\nonumber\ba
&((b_1-b_2)^2 + (b_3-b_4)^2 + (b_5-b_6)^2 + (b_7-b_8)^2)\\
&\cdot((b_1-b_2)^2 - (b_3-b_4)^2 - (b_5-b_6)^2 + (b_7-b_8)^2)\\
&\cdot(b_1 + b_2 + b_3 + b_4 + b_5 + b_6 + b_7 + b_8)\\
&\cdot(b_1 + b_2 + b_3 + b_4 - b_5 - b_6 - b_7 - b_8)\\
&\cdot(b_1 + b_2 - b_3 - b_4 + b_5 + b_6 - b_7 - b_8)\\
&\cdot(b_1 + b_2 - b_3 - b_4 - b_5 - b_6 + b_7 + b_8).
\ea\end{equation*}
\ignore{
Sage:

A4 = matrix(SR,8,8,[b1,b2,b3,b4,b5,b6,b7,b8,
        b2,b1,b4,b3,b6,b5,b8,b7,
        b4,b3,b1,b2,b7,b8,b6,b5,
        b3,b4,b2,b1,b8,b7,b5,b6,
        b6,b5,b8,b7,b1,b2,b3,b4,
        b5,b6,b7,b8,b2,b1,b4,b3,
        b8,b7,b6,b5,b3,b4,b1,b2,
        b7,b8,b5,b6,b4,b3,b2,b1]);

A4; A4.determinant().expand().factor();

Result:
(b1^2 - 2*b1*b2 + b2^2 + b3^2 - 2*b3*b4 + b4^2 + b5^2 - 2*b5*b6 + b6^2 + b7^2 - 2*b7*b8 + b8^2)*(b1^2 - 2*b1*b2 + b2^2 - b3^2 + 2*b3*b4 - b4^2 - b5^2 + 2*b5*b6 - b6^2 + b7^2 - 2*b7*b8 + b8^2)*(b1 + b2 + b3 + b4 + b5 + b6 + b7 + b8)*(b1 + b2 + b3 + b4 - b5 - b6 - b7 - b8)*(b1 + b2 - b3 - b4 + b5 + b6 - b7 - b8)*(b1 + b2 - b3 - b4 - b5 - b6 + b7 + b8)

}%end ignore

\subsection{}
It should be clear by now that the examples given can provide 
entertaining and instructive work and practice in linear algebra
for students ranging from high school up to upper division undergraduate
level.  The interested reader may be motivated to look for
other examples of finite group actions on the real line or the complex
plane, and to explore the corresponding functional equations.

It has to be conceded that most of the questions one can invent using these
methods will take too long to work out by hand to be useful in
competitions.  
But such questions can always be used as exercises for
college-level courses in the use of mathematical software such as
Sage, Maple or Mathematica. I used the open-source online Sagemath cell 
\cite{sagemath} 
to evaluate and factor some of the determinants above.

The next, short, section of this paper is aimed at readers with a more
advanced background in abstract algebra, and serves to round out
our story by locating it in relation to the history and theory of
group representations.

\section{The matrix and form}
When the group $H$ is trivial, the coefficients $a_{g,h}$
in Equation \eqref{E:general}
all take the form $a_{g,1}$, with $g\in G$, 
and we abbreviate this to $a_g$.
If, then, we index the coefficient matrix
of the system of equations \eqref{E:general-system}
by the elements of the finite group $G$, then the element in
the row indexed by $s$ and column indexed by $t$ 
is $a_{s^{-1}t}$.
Its determinant $\Delta$
 is independent of the order in which the elements
are listed.  The transposed matrix, with $(s,t)$ entry $a_{st^{-1}}$,
has the same determinant.
Group theorists call this determinant the
\emph{group determinant} of $G$. See
the historical account in
\cite{Curtis} and the well-named paper \cite{FS}.

For general finite groups $G$ and $H$, the form is
more complicated. 

A \emph{representation} of a group $G$ over some field $F$
is a homomorphism from
$G$ into the group $\GL(n,F)$ of invertible $n\times n$ 
matrices over $F$ (see, for instance,
\cite[Chapter 1]{Serre}, \cite{Hill}, or \cite{Bartel}).  Alternatively,
a representation may be thought of as a homomorphism from
$G$ into a group $\GL(V,F)$ of invertible linear endomorphisms
of a vector space $V$ over $F$. Representations are the single
most useful tool in advanced group theory. 
There is a form of degree $|G|$ associated to every
representation $\phi:G\to\GL(n,F)$.  It is defined by
$$  d(\phi)(\lambda):=
\det \phi\left(\sum_{g\in G} \lambda_g g \right). $$

The \emph{regular representation} of a finite group $G$ over $F$
is obtained by considering the elements of $G$ as
a basis for the $|G|$-dimensional
vector space $V$ of all formal sums
$$ \sum_{h\in G} \alpha_h h, $$
with each $\alpha_h\in F$, and mapping an element $g\in G$
to the matrix representing the invertible map
$$ \sum_{h\in G} \alpha_h h \mapsto 
\sum_{h\in G} \alpha_h gh. 
$$
 
The group determinant is the form associated to the
regular representation.
The factors of the group determinant over $\C$ are
the forms $d(\phi)$ corresponding to the irreducible
representations of $G$ over $\C$.

The group determinant is nonconstant, so vanishes only on
a proper subvariety of $F^G$ (the vector space of all functions from
$G$ into $F$, considered as an algebraic variety
\cite[Chapter 5]{Rotman}).  Whenever $F$ is infinite,
it is nonzero generically, and the equation usually has unique solution.

\section{Integral forms}
We conclude by noting that the same form gives another kind of problem, 
of which the following is an example.

\begin{question} Show that if two numbers  take the form
        $a^3+b^3+c^3-3abc$, where $a,b,c$ are integers, then their
        product has the same form.
\end{question}

(This once popped up on an intervarsity competition paper set by
\afriend.)
The question was easy for anyone who knew the factorization 
$$ a^3+b^3+c^3-3abc = (a+b+c)(a+b\omega+c\omega^2)(a+b\omega^2+c\omega),
$$
(--- including anyone who took the Irish Leaving Certificate
before the mid-sixties ---)
but is even easier when the form is recognized as a group determinant.
In fact, for any finite group $G$, with determinant $d$ over a field $F$,
we could consider the values taken by $d$ on $S^G$ (the set of all
functions from $G$ into $S$), where
$S$ is any subsemiring of $F$, i.e. a subset closed under addition
and multiplication. If we
denote this set of values by $d(S^G)$, then
$$ d(S^G) = \det \pi( SG ), $$
where $SG$ denotes the subset of the group ring $FG$ consisting
of elements with coefficients in $S$, and $\pi$ is
the regular representation.  Since $SG$ forms a subsemiring
of  $FG$, and $\det$ and $\pi$ are multiplicative,
it follows that $d(S^G)$ is closed under products.
For example, in characteristic zero one could insist
on positive integral coefficients, or coefficients
divisibly by $3$.

%Less abrupt closing paragraph:
Evidently, a similar problem may be posed about the group determinant
of any finite group $G$. For instance, $G=C_2\times C_2$
gives the form
\ignore{
	\begin{align*}
		&(a_1 1 + b_1\sigma + c_1\tau + \d_1\sigma\tau) 
(a_2 2 + b_2\sigma + c_2\tau + \d_2\sigma\tau)\\
		&(a_1a_2 + b_1 b_2 + c_1 c_2 + d_1d_2)1
		\\
		&(a_1b_2+a_2b_1 + c_1d_2+c_2d_1)\sigma
		\\
		&(a_1c_1+a_2c_1+b_1d_2+b_2d_1)\tau
		\\
		&(a_1d_2+a_2d_1+b_1c_2+b_2c_1)\sigma\tau
	\end{align*}
}%end ignore
$$ 
((a+b)^2-(c+d)^2)((a-b)^2-(c-d)^2). $$
Thus the set of values of this form, as $a,b,c,d$ range over all
positive integers, is closed under multiplication.
This kind of problem can always be tackled without sophistication,
and solved by pure ingenuity, which can be entertaining to
watch.

Incidentally, the abelian group $C_2\times C_2$ is isomorphic to that generated by the functions $-t$ and $1/t$, and the main method of this paper
gives the rather topical puzzle:
\begin{question} Suppose $f:\R\to\R$ and for each nonzero $x$
        satisfies the identity
        $$ f(x)+ 2 f(-x) + 4 f(1/x) + 8 f(-1/x) = 2025 x^2.
$$
        What is $f(3)$?
\end{question}
The reader is invited to verify that the answer is $-385$.
%(This nice round number comes about mainly because this
%group determinant, evaluated at $1,2,4,8$, equals $2025.)
\ignore{
	General solution with rhs 2025 x^2 is
	-45 \left(\frac{x^4-4}{x^2}\right)
}%end ignore

\subsection{acknowledgment}
I am grateful to \acolleague for useful conversation about the
content of this paper and to the editors and referee
for useful advice on the exposition.

\bibliographystyle{plain}
\bibliography{fnq}

\end{document}